\documentclass[12pt]{article}

\usepackage[T1]{fontenc}
\usepackage{avant}
\usepackage[cp1250]{inputenc}
\usepackage{amsmath,amssymb}
\usepackage{amsthm}
\usepackage{caption}
\usepackage{geometry}
\usepackage{float}

\usepackage{makeidx}
\usepackage[matrix,arrow,curve]{xy}
\usepackage{boxedminipage}
\usepackage{epic}
\usepackage{longtable}
\usepackage{ifthen}
\usepackage{tabularx}
\usepackage{parskip}
\usepackage{verbatim}
\usepackage{multirow}
\usepackage{graphicx}

\theoremstyle{plain}
\newtheorem{theorem}{Theorem}[section]
\newtheorem{proposition}[theorem]{Proposition}
\newtheorem{lemma}[theorem]{Lemma}
\newtheorem{definition}[theorem]{Definition}

\newtheorem{corollary}[theorem]{Corollary}

\usepackage{authblk}

\title{Proper $SL(2,\mathbb{R})$-actions on homogeneous spaces}

\author{Maciej Boche\'nski, Piotr Jastrz\k{e}bski, Takayuki Okuda and Aleksy Tralle}

\begin{document}

\maketitle{}

\newcommand{\leftrightdef}{\mathrel{\stackrel{\mathrm{def}}{\longleftrightarrow}}}
\newcommand{\simrightarrow}{\mathrel{\stackrel{\sim}{\rightarrow}}}
\newcommand{\simleftarrow}{\mathrel{\stackrel{\sim}{\leftarrow}}}
\newcommand{\bijarrow}{\mathrel{\stackrel{1:1}{\longleftrightarrow}}}
\newcommand{\ad}{\mathop{\mathrm{ad}}\nolimits}
\newcommand{\Ad}{\mathop{\mathrm{Ad}}\nolimits}
\newcommand{\ev}{\mathop{\mathrm{ev}}\nolimits}
\newcommand{\Exp}{\mathop{\mathrm{Exp}}\nolimits}
\newcommand{\Int}{\mathop{\mathrm{Int}}\nolimits}
\newcommand{\Aut}{\mathop{\mathrm{Aut}}\nolimits}
\newcommand{\End}{\mathop{\mathrm{End}}\nolimits}
\newcommand{\Map}{\mathop{\mathrm{Map}}\nolimits}
\newcommand{\Tr}{\mathop{\mathrm{Tr}}\nolimits}
\newcommand{\diag}{\mathop{\mathrm{diag}}\nolimits}
\newcommand{\level}{\mathrel{\mathrm{level}}}
\newcommand{\Ker}{\mathop{\mathrm{Ker}}\nolimits}
\newcommand{\locisom}{\mathrel{\stacrel{\mathrm{loc}}{=}}}
\newcommand{\Hom}{\mathop{\mathrm{Hom}}\nolimits}
\newcommand{\Lie}{\mathop{\mathrm{Lie}}\nolimits}
\newcommand{\rank}{\mathop{\mathrm{rank}}\nolimits}
\newcommand{\R}{\mathbb{R}}
\newcommand{\C}{\mathbb{C}}
\newcommand{\Q}{\mathbb{Q}}
\newcommand{\Z}{\mathbb{Z}}
\newcommand{\N}{\mathbb{N}}
\begin{abstract}
We study the existence problem of proper actions of $SL(2,\mathbb{R})$ on homogeneous spaces $G/H$ of reductive type. 
Based on Kobayashi's properness criterion [Math.~Ann.~(1989)], 
we show that $G/H$ admits a proper $SL(2,\mathbb{R})$-action via $G$ if a maximally split abelian subspace of Lie $H$ is included in the wall defined by a restricted root of Lie $G$. 
We also give a number of examples of such $G/H$.


\end{abstract}

\section{Introduction and main results}\label{sec:intro}

In the present paper, 
we find a number of examples of homogeneous spaces $G/H$ of reductive type which admit proper actions of $SL(2,\mathbb{R})$.
Our method is based on the criterion of proper actions established by Kobayashi \cite{kob2} (see Theorem \ref{tkob} below).

To begin with, let us recall the definition of proper action. 
Let $L$ be a locally compact topological group acting continuously on a locally compact Hausdorff topological space $X$. 
This action is said to be {\it proper} if for every compact subset $S \subset X$ the set
$$L_{S}:=\{  g\in L \ | \ g\cdot S \cap S \neq \emptyset \}$$
is compact. 
The action is called {\it properly 
discontinuous} if it is proper and the topology of $L$ is discrete.

A general study of discontinuous groups for homogeneous spaces $G/H$ with $H$ non-compact was initiated by Kobayashi \cite{kob2} in full generality. 
Following the survey \cite{kobbb2}, 
we recall some of the general problems and perspectives as below. 
Let $G$ be a Lie group and $H$ a closed subgroup of $G$.
A discrete subgroup $\Gamma$ of $G$ is said to be a {\it discontinuous group } for the homogeneous space $G/H$ if the natural $\Gamma$-action on $G/H$ is properly discontinuous and fixed-point free.
For such $\Gamma$, the quotient space $\Gamma \backslash G/H$ is a smooth manifold to be referred to as a Clifford--Klein form (see \cite{kob2}). 
In more detail $\Gamma \backslash G/H$ has a unique smooth manifold structure such that the quotient map $G/H \rightarrow \Gamma \backslash G/H$ is a $C^{\infty}$-covering. 
Thus any $G$-invariant local geometric structure on $G/H$ descends onto $\Gamma \backslash G/H$. 
This opens a perspective of applying results on the existence of such groups in geometry.
 If $\Gamma \backslash G/H$ is compact, we say that $G/H$ admits a compact Clifford--Klein form. 

The problem of finding discontinuous groups for $G/H$ is straightforward when $H$ is compact. 
In such cases any discrete subgroup of $G$ acts properly discontinuously on $G/H$. 
So the interesting case is when $H$ is non-compact. 
For instance, if $\text{rank}_{\mathbb{R}}G=\text{rank}_{\mathbb{R}}H$ then only finite groups can act properly discontinuously on $G/H$ (it is the Calabi-Markus phenomenon, see \cite{cam} and \cite{kob2}). 

One of effective ways of constructing properly discontinuous groups for $G/H$ is to find closed Lie subgroups $L$ of $G$ acting properly on $G/H$.
Then any torsion-free discrete subgroup $\Gamma$ of $L$ is a discontinuous group for $G/H$.
For homogeneous spaces $G/H$ of reductive type, 
many results for closed Lie subgroups $L$ of $G$ acting properly on $G/H$ 
can be found in \cite{ben}, \cite{kass}, \cite{kass-kob}, \cite{kobaaa}, \cite{kobbb}, \cite{ko}, \cite{rs}, \cite{ok}, \cite{d}.
One of the most important results in this area is Kobayashi's criterion \cite{kob2} of the properness of the $L$-action on $G/H$ in terms of Lie algebras of $G$, $H$ and $L$ (see Section \ref{sec:triples} for more details).

We are interested in homogeneous spaces $G/H$ of reductive type which admits proper actions of $SL(2,\mathbb{R})$.
Throughout this paper, we use the following terminology:

\begin{definition}\label{def:properviaG}
We say that 
a homogeneous space $G/H$ admits a proper $SL(2,\R)$-action via $G$
if there exists a Lie group homomorphism $\rho : SL(2,\R) \rightarrow G$ such that the action of $SL(2,\R)$ on $G/H$ induced by $\rho$ is proper. 
\end{definition}

In Definition \ref{def:properviaG}, 
the image $\rho(SL(2,\R))$ in $G$ of $\rho$ is isomorphic to $SL(2,\R)$ or $PSL(2,\R)$.
Note that by the lifting theorem of surface groups (cf.~\cite{kra}), 
discontinuous groups $\Gamma$ for $G/H$ which is isomorphic to the surface group $\pi_1(\Sigma_g)$ for some $g \geq 2$
can be constructed as discrete subgroups of $\rho(SL(2,\R))$.

Kulkarni \cite{rs} studies parameters $(p,q)$ such that $S^{p,q} = O(p+1,q)/O(p,q)$ admitting proper $SL(2,\R)$-actions via $O(p+1,q)$.
In general, semisimple symmetric spaces $G/H$ admitting proper $SL(2,\mathbb{R})$-actions via $G$ are classified by Teduka \cite{d} and the third author \cite{ok}.
Continuing this line of research, in \cite{bt}, 
the first and the fourth authors give a sufficient condition on the Lie algebra data of the pair $(\mathfrak{g},\mathfrak{h})$ which implies the existence of 
non virtually-abelian discontinuous groups for $G/H$. 
This enabled them to extend results of \cite{ben} to homogeneous spaces generated by automorphisms of order $3$.

In this paper, based on Kobayashi's properness criterion \cite{kob2}, we show a sufficient condition on homogeneous spaces $G/H$ for the existence of proper $SL(2,\mathbb{R})$-actions via $G$. 
Our main theorems are stated in Theorems \ref{twggg} and \ref{twggg1} below.  

Let $G$ be a connected non-compact real semisimple linear Lie group,
$\theta$ a Cartan involution on $G$ 
and $H$ a $\theta$-stable closed subgroup of $G$.
Then $G/H$ is of reductive type.
We take maximally split abelian subspaces $\mathfrak{a}_{\mathfrak{g}}$, $\mathfrak{a}_{\mathfrak{h}}$ of $\mathfrak{g} := \textrm{Lie} \ G$, $\mathfrak{h} := \textrm{Lie} \ H$ respectively, with $\mathfrak{a}_{\mathfrak{h}} \subset \mathfrak{a}_{\mathfrak{g}}$ (see Section \ref{sec:triples} for the notation).
The restricted root system of $(\mathfrak{g},\mathfrak{a}_{\mathfrak{g}})$ is denoted by $\Sigma_{\mathfrak{g}}\subset \mathfrak{a}_{\mathfrak{g}}^*$. 

\begin{theorem}\label{twggg}
Assume that there exists $\alpha \in \Sigma_{\mathfrak{g}}$ such that $\alpha (\mathfrak{a}_{\mathfrak{h}})= \{ 0 \}$. 
Then the homogeneous space $G/H$ admits a proper $SL(2,\mathbb{R})$-action via $G$.
\end{theorem}

We remark that the assumption of Theorem \ref{twggg} depends only on $(\mathfrak{g},\mathfrak{h})$ but not on the choice of $\mathfrak{a}_{\mathfrak{g}}$ and $\mathfrak{a}_{\mathfrak{h}} \subset \mathfrak{a}_{\mathfrak{g}}$.

The next result yields a class of homogeneous spaces for which the condition of Theorem \ref{twggg} is easy to check.

\begin{theorem}\label{twggg1}
For each $\alpha \in \Sigma_{\mathfrak{g}}$, 
we take a closed Lie subgroup $H = H^{\alpha}$ of $G$ such that $\textrm{Lie} \ H = [Z_{\mathfrak{g}}(A_{\alpha}),Z_{\mathfrak{g}}(A_{\alpha})]$, where $A_{\alpha}$ is the coroot of $\alpha$. 
Then the homogeneous space $G/H$ admits a proper $SL(2,\mathbb{R})$-action via $G$.
\end{theorem}

Note that in the setting of Theorem \ref{twggg1}, $\mathfrak{h} = \textrm{Lie} \ H$ is the semisimple part of the Levi subalgebra $Z_{\mathfrak{g}}(A_{\alpha})$ of $\mathfrak{g}$.

Families of examples of Theorem \ref{twggg} and \ref{twggg1} are given in Section \ref{section:ex}.
In greater detail, let $\mathfrak{g}$ be a split real form of a complex semisimple Lie algebra and $\mathfrak{h}$ a regular semisimple subalgebra of $\mathfrak{g} := \Lie G$.
Using the classification results of Oshima \cite{os} 
on $S$-closed subsystems of root systems, 
we can always check whether or not $(\mathfrak{g},\mathfrak{h})$ satisfies the assumption of Theorem \ref{twggg}.
Then, as in Tables \ref{tabCl} and \ref{tabEx} below, 
we obtain a number of examples $(\mathfrak{g},\mathfrak{h})$ satisfying the assumption in Theorem \ref{twggg}.
Furthermore, for non-split $\mathfrak{g}$, we also obtain some examples of $(\mathfrak{g},\mathfrak{h})$ in Theorem \ref{twggg1} as in Table \ref{tab1} below.

We also study homogeneous spaces $G/H$ of reductive type admitting proper $SL(2,\mathbb{R})$-actions via $G$ in the cases where the center of $H$ is compact and $\text{rank}_{\mathbb{R}}(\mathfrak{h}) = 1$ (see Section \ref{subsec:rank-one}).

\vskip6pt
\noindent {\bf Acknowledgment}. 
The third author is supported by JSPS KAKENHI Grant Number JP16K17594 in this work.


\section{Kobayashi's properness criterion}\label{sec:triples}

In this section, we recall Kobayashi's properness criterion \cite[Theorem 4.1]{kob2} in a form that we shall need.

Let $G$ be a connected real semisimple Lie group with its Lie algebra $\mathfrak{g}$. 
We fix a Cartan decomposition $\mathfrak{g} = \mathfrak{k} + \mathfrak{p}$ and a maximal abelian subspace $\mathfrak{a}_\mathfrak{g}$ of $\mathfrak{p}$.
Such subspace $\mathfrak{a}_\mathfrak{g}$ of $\mathfrak{g}$ is unique up to the adjoint action of $G$ and called a {\it maximally split abelian subspace} of $\mathfrak{g}$.
Let us denote by $\Sigma_\mathfrak{g}$ and $W_{\mathfrak{g}}$ the restricted root system of $(\mathfrak{g},\mathfrak{a}_\mathfrak{g})$ and the Weyl group of $\Sigma_\mathfrak{g}$ acting on $\mathfrak{a}_\mathfrak{g}$, respectively.

Throughout this paper,
a closed Lie subgroup $H$ of $G$ is said to be {\it reductive in} $G$ 
if there exists a Cartan involution $\theta$ on $G$ such that $H$ is $\theta$-stable.
Then $H$ is reductive as a Lie group (the converse is not true in general).
For such a closed subgroup $H$ of $G$, 
$\mathfrak{a}_\mathfrak{h}$ denotes a subspace of $\mathfrak{a}_\mathfrak{g}$ which is conjugate to a maximally split abelian subspace of $\mathfrak{h}$.
That is, there exists $g \in G$ and a Cartan involution $\theta$ on $G$ with $\theta(H) = H$ such that $\Ad(g) \mathfrak{a}_\mathfrak{h} \subset \mathfrak{p}_\mathfrak{h}$ where $\mathfrak{h} = \mathfrak{k}_\mathfrak{h} + \mathfrak{p}_\mathfrak{h}$ denotes the Cartan decomposition of $\mathfrak{h} = \Lie H$ by $\theta$.
Such a subspace $\mathfrak{a}_\mathfrak{h}$ of $\mathfrak{a}_\mathfrak{g}$ is unique up to the conjugations by the Weyl gourp $W_\mathfrak{g}$.

Let us fix closed subgroups $H$ and $L$ which are both reductive in $G$.
The Lie algebras of $H$ and $L$ are denoted by $\mathfrak{h}$ and $\mathfrak{l}$, respectively.
We fix subspaces $\mathfrak{a}_\mathfrak{h}$ and $\mathfrak{a}_\mathfrak{l}$ of $\mathfrak{a}_\mathfrak{g}$ as above.

The following simple criterion by Kobayashi \cite{kob2} for the properness of the action of $L$ on $G/H$ is the crucial tool of this paper.

\begin{theorem}[{\cite[Theorem 4.1]{kob2}}]\label{tkob}
The following three conditions on $(G,H,L)$ are equivalent:
\begin{description}
\item[(i)] The $L$-action on $G/H$ is proper.
\item[(ii)] The $H$-action on $G/L$ is proper.
\item[(iii)]  $\mathfrak{a}_{\mathfrak{h}} \cap W_{\mathfrak{g}} \mathfrak{a}_{\mathfrak{l}} = \{ 0 \}$ in $\mathfrak{a}_\mathfrak{g}$.
\end{description}
\end{theorem}

Let us consider the cases where $\mathfrak{l} \simeq \mathfrak{sl}(2,\R)$.
A triple $(A,X,Y)$ of vectors in $\mathfrak{g}$ is called the \textbf{{\it $\mathfrak{sl}_2$-triple}} if
$[A,X]=2X$, $[A,Y]=-2Y$ and $[X,Y]=A$.
Then there exists a Lie group homomorphism 
$\Phi : SL(2,\mathbb{R}) \rightarrow G$
with \[
d\Phi  \left( \begin{array}{cc} 1 & 0  \\ 0 & -1  \end{array} \right) = A, \ \ d\Phi  
\left( \begin{array}{cc} 0 & 1  \\ 0 & 0  \end{array} \right) = X, \ \ d\Phi  \left( 
\begin{array}{cc} 0 & 0  \\ 1 & 0  \end{array} \right) = Y.
\]
Indeed, the Lie algebra homomorphism $d\Phi : \mathfrak{sl}(2,\R) \rightarrow \mathfrak{g}$ defined by the $\mathfrak{sl}_2$-triple $(A,X,Y)$ can be lifted to $\Phi: SL(2,\mathbb{R}) \rightarrow G$ since the complexification $SL(2,\mathbb{C})$ of $SL(2,\mathbb{R})$ is simply-connected.

As a corollary to Theorem \ref{tkob} for $\mathfrak{l} \simeq \mathfrak{sl}(2,\R)$, 
we have the following claim:

\begin{corollary}[{\cite[Corollary 5.5]{ok}}]\label{tthh}
Let $G$ be a connected semisimple linear Lie group and $H$ a closed subgroup of $G$ which is redctive in $G$. 
For a given $\mathfrak{sl}_2$-triple $(A,X,Y)$ in $\mathfrak{g}$, 
we denote by $\rho : SL(2,\mathbb{R}) \rightarrow G$ the corresponding Lie group homomorphism. 
Then the following two conditions on $(G,H,\rho)$ are equivalent:
\begin{description}
	\item[(i)] The $SL(2,\mathbb{R})$-action on $G/H$ induced by $\rho$ is proper.
	\item[(ii)] The adjoint orbit through $A$ in $\mathfrak{g}$ does not meet $\mathfrak{h}$.
\end{description}
\end{corollary}

It should be remarked that in the case where $A \in \mathfrak{a}$, the second condition in Theorem \ref{tthh} is equivalent to the condition that $A\notin W\cdot \mathfrak{a}_{\mathfrak{h}}$.

\section{Proofs of Theorem \ref{twggg} and Theorem \ref{twggg1}}\label{sec:nilpotent}
\label{sec:onil}

Let $\mathfrak{g}$ be a non-compact real semisimple Lie algebra. We take a maximally split abelian subalgebra $\mathfrak{a}_{\mathfrak{g}}$ and denote by $\Sigma_{\mathfrak{g}} \subset \mathfrak{a}_{\mathfrak{g}}^*$ the restricted root system of $(\mathfrak{g},\mathfrak{a}_{\mathfrak{g}})$ as before.

We show the following lemma:

\begin{lemma}\label{lemma:principal}
There exists an $\mathfrak{sl}_2$-triple $(A,X,Y)$ in $\mathfrak{g}$ such that  $\alpha (A) \neq 0$ for any $\alpha \in \Sigma_{\mathfrak{g}}$.
\end{lemma}

It should be noted that 
the claim of the lemma above can be found in \cite[Section 4]{Djo} without proofs.
For the convenience of the readers, we give a proof as follows.

\begin{proof}[Proof of Lemma \ref{lemma:principal}]
The proof is similar to that of in the complex case (see \cite[Theroem 4.1.6]{cmg}). Let us fix a simple system $\Pi_{\mathfrak{g}}$ of the restricted root system $\Sigma_{\mathfrak{g}}.$ We denote by $A \in \mathfrak{a}_{\mathfrak{g}}$ the unique element in $\mathfrak{a}_{\mathfrak{g}}$ such that $\alpha (A)=2$ for all $\alpha \in \Pi_{\mathfrak{g}}.$ It is enough to show that there exists $X,Y\in \mathfrak{g}$ such that $(A,X,Y)$ is an $\mathfrak{sl}_2$-triple in $\mathfrak{g}$.

For each $\alpha \in \Pi_{\mathfrak{g}}$ we denote by $A_{\alpha}$ the coroot of $\alpha$. 
Then it is well-known that there exist $X_{\alpha} \in \mathfrak{g}_{\alpha}$ and $Y_{\alpha} \in \mathfrak{g}_{-\alpha}$ such that $(A_{\alpha}, X_{\alpha}, Y_{\alpha})$ is an $\mathfrak{sl}_2$-triple in $\mathfrak{g}$. 
Here $\mathfrak{g}_{\pm \alpha}$ is the root space of $\pm \alpha$ in $\mathfrak{g}$ (not in $\mathfrak{g}^{\mathbb{C}}$). Since $\{ A_{\alpha} \ | \ \alpha \in \Pi_{\mathfrak{g}} \}$ is a basis of $\mathfrak{a}_{\mathfrak{g}}$ the element $A\in \mathfrak{a}_{\mathfrak{g}}$ can be written as $A = \sum_{\alpha \in \Pi_{\mathfrak{g}}} a_{\alpha} A_{\alpha}$ for some $a_{\alpha}\in \mathbb{R}$ $(\alpha \in \Pi_{\mathfrak{g}})$. Let us put
$$X:=\sum_{\alpha \in \Pi_{\mathfrak{g}}}X_{\alpha} \ \textrm{and} \ Y:= \sum_{\alpha \in \Pi_{\mathfrak{g}}} a_{\alpha}Y_{\alpha}.$$
We shall prove that $(A,X,Y)$ is an $\mathfrak{sl}_2$-triple.
We see that 
$$[A,X] = \sum\limits_{\alpha \in \Pi_{\mathfrak{g}}} \alpha (A) X_{\alpha}=2X,$$
$$[A,Y] = \sum\limits_{\alpha \in \Pi_{\mathfrak{g}}} a_{\alpha}(-\alpha (A)) Y_{\alpha}=-2Y.$$
Recall that since $\Pi_{\mathfrak{g}}$ is a simple system, $\alpha -\beta$ is not a root in $\Sigma_{\mathfrak{g}}$ for any $\alpha ,\beta \in \Sigma_{\mathfrak{g}}$. Therefore $[X_{\alpha}, Y_{\beta}]=0$ if $\alpha\neq \beta$ and hence we have
$$[X,Y] = \sum\limits_{\alpha, \beta \in \Pi_{\mathfrak{g}}} a_{\beta}[X_{\alpha},Y_{\beta}]=\sum\limits_{\alpha \in \Pi_{\mathfrak{g}}} a_{\alpha}[X_{\alpha},Y_{\alpha}]=\sum\limits_{\alpha \in \Pi_{\mathfrak{g}}} a_{\alpha} A_{\alpha}=A. $$
This completes the proof.
\end{proof}

Now we are ready to prove Theorem \ref{twggg}.

\begin{proof}[Proof of Theorem \ref{twggg}]
There exists $\alpha \in \Sigma_{\mathfrak{g}}$ such that $\alpha (\mathfrak{a}_{\mathfrak{h}}) =\{ 0 \}$ and hence $\mathfrak{a}_{\mathfrak{h}}$ is contained in the wall $C_{\alpha} := \{ V\in \mathfrak{a}_{\mathfrak{g}} \ | \ \alpha(V)=0 \}$ defined by $\alpha$. 
By Lemma \ref{lemma:principal}, there exists an $\mathfrak{sl} (2,\mathbb{R})$-triple $(A,X,Y)$ in $\mathfrak{g}$ such that $\beta (A)\neq 0$ for any $\beta \in \Sigma_{\mathfrak{g}}$. Therefore $A$ does not lie in the walls in $\mathfrak{a}_{\mathfrak{g}}$. Hence we see that the orbit of the Weyl group through $A$ does not meet $\mathfrak{a}_{\mathfrak{h}}$. The proof is completed by Corollary \ref{tthh}.
\end{proof}

Using Theorem \ref{twggg}, we prove Theorem \ref{twggg1} as follows:

\begin{proof}[Proof of Theorem \ref{twggg1}]
Recall that $A_{\alpha}$ is a coroot of $\alpha \in \Sigma_{\mathfrak{g}}$ and $\mathfrak{h} = \textrm{Lie} \ H = [Z_{\mathfrak{g}}(A_{\alpha}),Z_{\mathfrak{g}}(A_{\alpha})]$. 
Therefore $\{ V \in \mathfrak{a}_{\mathfrak{g}} \ | \ \alpha (V)=0 \}$ is a maximally split abelian subspace of $\mathfrak{h}.$ Thus Theorem \ref{twggg1} follows from Theorem \ref{twggg}.
\end{proof}

\section{Examples}\label{section:ex}

\subsection{Examples from Theorem \ref{twggg}, the split case}\label{sub51}

Let $\Sigma$ be a reduced root system.
As in \cite{os}, 
we say that a subset $\Xi$ of $\Sigma$ is an {\it $S$-closed subsystem}
if $\alpha + \beta \in \Xi$ for any $\alpha, \beta \in \Xi$ with $\alpha + \beta \in \Sigma$.
Note that any $S$-closed subsystem $\Xi$ of $\Sigma$ is a subsystem of $\Sigma$ in the sense of \cite{os}, that is, $s_\alpha(\Xi) = \Xi$ for any $\alpha \in \Xi$, where $s_\alpha$ is the reflection on $\mathfrak{a}^*$ through the hyperplane perpendicular to $\alpha$.

Let us fix a reduced root system $\Sigma$ and an $S$-closed subsystem $\Xi$ of $\Sigma$.
We denote by $\mathfrak{g}_\Sigma$ a split real form of a complex semisimple Lie algebra corresponding to $\Sigma$.
We fix a maximally split abelian subspace $\mathfrak{a}$ of $\mathfrak{g}_\Sigma$
and then $\Sigma$ is realized as the restricted root system of $(\mathfrak{g}_\Sigma,\mathfrak{a})$.
The root space decomposition of $\mathfrak{g}_\Sigma$ by $\mathfrak{a}$ can be written as 
\[
\mathfrak{g}_\Sigma = \mathfrak{a} + \sum_{\lambda \in \Sigma} (\mathfrak{g}_\Sigma)_{\lambda}.
\] 
For the $S$-closed subsystem $\Xi$ of $\Sigma$, 
we define a semisimple subalgebra $\mathfrak{h}_\Xi$ by
\[
\mathfrak{h}_{\Xi} := \mathfrak{a}_{\Xi} + \sum_{\alpha \in \Xi} (\mathfrak{g}_\Sigma)_{\alpha},
\]
where $\mathfrak{a}_\Xi$ is the subspace of $\mathfrak{a}$ spaned by the set of all coroots of $\alpha \in \Xi$.
Then $\mathfrak{h}_\Xi$ is an $\mathfrak{a}$-regular semisimple subalgebra of $\mathfrak{g}_\Sigma$ in the sence of \cite{DFG}, that is, $\mathfrak{h}_\Xi$ is semisimple and normalized by $\mathfrak{a}$.
Conversely, any $\mathfrak{a}$-regular semisimple subalgebra of $\mathfrak{g}_\Sigma$ can be obtained as $\mathfrak{h}_\Xi$ for an $S$-closed subsystem $\Xi$ of $\Sigma$.

Let us take a connected semisimple linear Lie group $G_\Sigma$ and its closed subgroup $H_\Xi$ with $\Lie G_\Sigma = \mathfrak{g}_\Sigma$ and $\Lie H_\Xi = \mathfrak{h}_\Xi$.
Then by Theorem \ref{twggg}, 
$G_\Sigma/H_{\Xi}$ admits a proper $SL(2,\mathbb{R})$-action via $G_\Sigma$
if 
\[
\Xi^{\perp} := \{ \alpha \in \Sigma \ | \ \left \langle \alpha , \beta \right \rangle = 0 \ \textrm{for any} \ \beta \in \Xi \} \not=\emptyset.
\]

In \cite{os}, the classification of $S$-closed subsystems $\Xi$ of irreducible reduced root systems $\Sigma$ with $\Xi^\perp \neq \emptyset$ can be found in Tables 10.1 and 10.2. 
Using this classification, we obtain a number of examples 
of homogeneous spaces $G_\Sigma/H_\Xi$ admitting proper $SL(2,\mathbb{R})$-actions via $G_\Sigma$.
In Tables \ref{tabCl} and \ref{tabEx}, 
we give a complete list of pairs of irreducibe reduced root systems $(\Sigma,\Xi)$ such that $\Xi$ can be realized as an $S$-closed subsystem of $\Sigma$ with $\Xi^{\perp} \neq \emptyset$.

\begin{center}
 \begin{table}[h]
 \centering
 {\footnotesize
\begin{tabular}{|c|c|c|c|}
\hline $\Sigma$ & $\Xi$ & $\Xi^{\perp}$ & $\mathfrak{g}_\Sigma$, $\mathfrak{h}_\Xi$ \\ \hline
$A_{n-1}$ & $A_{k-1}$ &$A_{n-k-1}$ & $\mathfrak{sl}(n,\mathbb{R}), \ \mathfrak{sl}(k,\mathbb{R}), \ \  2\leq k\leq n-2$
\\ \hline
$D_n$ & $A_{k-1}$ &$D_{n-k}$ & $\mathfrak{so}(n,n), \ \mathfrak{sl}(k,\mathbb{R}), \ \ 2\leq k \leq n-2, \ \ n\geq 5$
\\ \hline
$D_n$ & $D_k$ &$D_{n-k}$ & $\mathfrak{so}(n,n), \ \mathfrak{so}(k,k), \ \ 4< k \leq n-2, \ \ n\geq 6$
\\ \hline
$B_n$ & $A_{k-1}$ &$B_{n-k}$ & $\mathfrak{so}(n,n+1), \ \mathfrak{sl}(k,\mathbb{R}), \ \ 2\leq k \leq n-1, \ \ n\geq 5$
\\ \hline
$B_n$ & $D_k$ &$B_{n-k}$ & $\mathfrak{so}(n,n+1), \ \mathfrak{so}(k,k), \ \ 4\leq k \leq n-1, \ \ n\geq 5$
\\ \hline
$B_n$ & $B_k$ &$B_{n-k}$ & $\mathfrak{so}(n,n+1), \ \mathfrak{so}(k,k+1), \ \ 2\leq k \leq n-1, \ \ n\geq 3$
\\ \hline
$C_n$ & $C_k$ &$C_{n-k}$ & $\mathfrak{sp}(n,\mathbb{R}), \ \mathfrak{sp}(k,\mathbb{R}), \ \ 1\leq k \leq n-1, \ \ n\geq 2$
\\ \hline
$C_n$ & $D_k$ &$C_{n-k}$ & $\mathfrak{sp}(n,\mathbb{R}), \ \mathfrak{so}(k,k), \ \ 3\leq k \leq n-1, \ \ n\geq 5$
\\ \hline
$C_n$ & $A_{k-1}$ &$C_{n-k}$ & $\mathfrak{sp}(n,\mathbb{R}), \ \mathfrak{sl}(k,\mathbb{R}), \ \ 1\leq k \leq n-2, \ \ n\geq 3$ \\ \hline
\end{tabular}
}
\captionsetup{justification=centering}
 \caption{
List of pairs of irreducibe reduced root systems $(\Sigma,\Xi)$ with $\Sigma$ of classical type
such that $\Xi$ can be realized as an $S$-closed subsystem of $\Sigma$ with $\Xi^{\perp} \neq \emptyset$.
 }
 \label{tabCl}
\end{table}
\end{center}

\begin{center}
 \begin{table}[h]
 \centering
 {\footnotesize
\begin{tabular}{|c|c|c|c|}
\hline $\Sigma$ & $\Xi$ & $\Xi^{\perp}$ & $\mathfrak{g}_\Sigma$, $\mathfrak{h}_\Xi$ \\ \hline
$E_6$ & $A_1$ &$A_5$ & $\mathfrak{e}_{6}^{\text{I}}, \ \mathfrak{sl}(2,\mathbb{R})$
\\ \hline
$E_6$ & $A_2$ &$2A_2$ & $\mathfrak{e}_{6}^{\text{I}}, \ \mathfrak{sl}(3,\mathbb{R})$
\\ \hline
$E_6$ & $A_3$ &$2A_1$ & $\mathfrak{e}_{6}^{\text{I}}, \ \mathfrak{sl}(4,\mathbb{R})$
\\ \hline
$E_6$ & $A_4$ &$A_1$ & $\mathfrak{e}_{6}^{\text{I}}, \ \mathfrak{sl}(5,\mathbb{R})$
\\ \hline
$E_6$ & $A_5$ &$A_1$ & $\mathfrak{e}_{6}^{\text{I}}, \ \mathfrak{sl}(6,\mathbb{R})$
\\ \hline
$E_7$ & $A_1$ &$D_6$ & $\mathfrak{e}_{7}^{\text{V}}, \ \mathfrak{sl}(2,\mathbb{R})$
\\ \hline
$E_7$ & $A_2$ &$A_5$ & $\mathfrak{e}_{7}^{\text{V}}, \ \mathfrak{sl}(3,\mathbb{R})$
\\ \hline
$E_7$ & $A_3$ &$A_3+A_1$ & $\mathfrak{e}_{7}^{\text{V}}, \ \mathfrak{sl}(4,\mathbb{R})$
\\ \hline
$E_7$ & $A_4$ &$A_2$ & $\mathfrak{e}_{7}^{\text{V}}, \ \mathfrak{sl}(5,\mathbb{R})$
\\ \hline
$E_7$ & $A_5$ &$A_2$ & $\mathfrak{e}_{7}^{\text{V}}, \ \mathfrak{sl}(6,\mathbb{R})$
\\ \hline
$E_7$ & $D_4$ &$3A_1$ & $\mathfrak{e}_{7}^{\text{V}}, \ \mathfrak{so}(4,4)$
\\ \hline
$E_7$ & $D_5$ &$A_1$ & $\mathfrak{e}_{7}^{\text{V}}, \ \mathfrak{so}(5,5)$
\\ \hline
$E_7$ & $D_6$ &$A_1$ & $\mathfrak{e}_{7}^{\text{V}}, \ \mathfrak{so}(6,6)$
\\ \hline
$E_8$ & $A_1$ &$E_7$ & $\mathfrak{e}_{8}^{\text{VIII}}, \ \mathfrak{sl}(2,\mathbb{R})$
\\ \hline
$E_8$ & $A_2$ &$E_6$ & $\mathfrak{e}_{8}^{\text{VIII}}, \ \mathfrak{sl}(3,\mathbb{R})$

\\\hline
\end{tabular}
\begin{tabular}{|c|c|c|c|}
\hline $\Sigma$ & $\Xi$ & $\Xi^{\perp}$ & $\mathfrak{g}_\Sigma$, $\mathfrak{h}_\Xi$ 
\\ \hline
$E_8$ & $A_3$ &$D_5$ & $\mathfrak{e}_{8}^{\text{VIII}}, \ \mathfrak{sl}(4,\mathbb{R})$
\\ \hline
$E_8$ & $A_4$ &$A_4$ & $\mathfrak{e}_{8}^{\text{VIII}}, \ \mathfrak{sl}(5,\mathbb{R})$
\\ \hline
$E_8$ & $A_5$ &$A_2+A_1$ & $\mathfrak{e}_{8}^{\text{VIII}}, \ \mathfrak{sl}(6,\mathbb{R})$
\\ \hline
$E_8$ & $A_6$ &$A_1$ & $\mathfrak{e}_{8}^{\text{VIII}}, \ \mathfrak{sl}(7,\mathbb{R})$
\\ \hline
$E_8$ & $A_7$ &$A_1$ & $\mathfrak{e}_{8}^{\text{VIII}}, \ \mathfrak{sl}(8,\mathbb{R})$
\\ \hline

$E_8$ & $D_4$ &$D_4$ & $\mathfrak{e}_{8}^{\text{VIII}}, \ \mathfrak{so}(4,4)$
\\ \hline
$E_8$ & $D_5$ &$A_3$ & $\mathfrak{e}_{8}^{\text{VIII}}, \ \mathfrak{so}(5,5)$
\\ \hline
$E_8$ & $D_6$ &$2A_1$ & $\mathfrak{e}_{8}^{\text{VIII}}, \ \mathfrak{so}(6,6)$
\\ \hline
$E_8$ & $E_6$ &$A_2$ & $\mathfrak{e}_{8}^{\text{VIII}}, \ \mathfrak{e}_{6}^{\text{I}}$
\\ \hline
$E_8$ & $E_7$ &$A_1$ & $\mathfrak{e}_{8}^{\text{VIII}}, \ \mathfrak{e}_{7}^{\text{V}}$
\\ \hline
$F_4$ & $A_1$ &$C_3$ & $\mathfrak{f}_{4}^{\text{I}}, \ \mathfrak{sl}(2,\mathbb{R})$
\\ \hline
$F_4$ & $A_2$ &$A_2$ & $\mathfrak{f}_{4}^{\text{I}}, \ \mathfrak{sl}(3,\mathbb{R})$
\\ \hline
$F_4$ & $B_2$ &$B_2$ & $\mathfrak{f}_{4}^{\text{I}}, \ \mathfrak{so}(2,3)$
\\ \hline
$F_4$ & $B_3$ &$A_1$ & $\mathfrak{f}_{4}^{\text{I}}, \ \mathfrak{so}(3,4)$
\\ \hline
$G_2$ & $A_1$ &$A_1$ & $\mathfrak{g}_{2}, \ \mathfrak{sl}(2,\mathbb{R})$ \\ \hline
\end{tabular}
 }
\captionsetup{justification=centering}
 \caption{
List of pairs of irreducibe reduced root systems $(\Sigma,\Xi)$ with $\Sigma$ of exceptional type
such that $\Xi$ can be realized as an $S$-closed subsystem of $\Sigma$ with $\Xi^{\perp} \neq \emptyset$.
 }
 \label{tabEx}
\end{table}
\end{center}

\subsection{Examples from Theorem \ref{twggg1}}\label{sub52}
We describe a procedure to construct examples of homogeneous spaces $G/H$ of reductive type admitting proper $SL(2,\mathbb{R})$-actions via $G$ by Theorem \ref{twggg1}. 
Every step will be illustrated for an example $\mathfrak{g}=\mathfrak{su}^{*}(10)$ as below.

\textbf{Step 1.} Draw the Satake diagram $S_{\mathfrak{g}}$ of $\mathfrak{g}$:

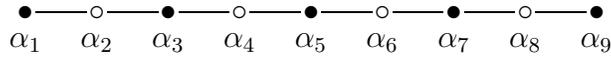
$$ \xymatrix@1@R=2pt@!C=3pt{
{\bullet}  \ar@{-}[r]& {\circ} \ar@{-}[r]&{\bullet}  \ar@{-}[r]& {\circ} \ar@{-}[r] &{\bullet}  \ar@{-}[r]& {\circ}\ar@{-}[r] &{\bullet}  \ar@{-}[r]& {\circ} \ar@{-}[r] &{\bullet}\\ \alpha_1 &\alpha_2 &\alpha_3&\alpha_4&\alpha_5&\alpha_6&\alpha_7&\alpha_8&\alpha_9
 }$$
Here $\{  \alpha_{1},...,\alpha_{9} \}$ is a simple system of the reduced root system $\Delta$ of the complexification $\mathfrak{sl}(10,\mathbb{C})$ of $\mathfrak{su}^{*}(10)$.

\textbf{Step 2.} Draw the extended Dynkin diagram $\mathrm{EDD}_{\Sigma}$ of the restricted root system $\Sigma_{\mathfrak{g}}$ of $\mathfrak{g}$:

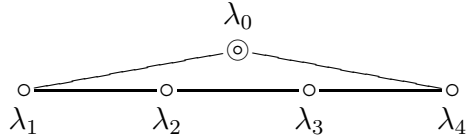
$$ \xymatrix@1@R=2pt@!C=3pt{ &&&\lambda_0 \\ &&&{\circledcirc} \ar@{-}[ddlll] \ar@{-}[ddrrr]\\ \\
{\circ}  \ar@{-}[rr] && {\circ} \ar@{-}[rr]&&{\circ}  \ar@{-}[rr]&& {\circ} \\ \lambda_1 &&\lambda_2 &&\lambda_3&&\lambda_4
 } $$
Here $\Pi_{\mathfrak{g}}:= \{  \lambda_{1}, \lambda_{2}, \lambda_{3}, \lambda_{4} \} $ is a simple system of $\Sigma_{\mathfrak{g}}$ with $r(\alpha_{2k})=\lambda_{k}$ for $k=1,2,3,4$ by the restriction map $r:\Delta \rightarrow \Sigma \cup \{ 0 \}, $ and $\lambda_{0}$ is the lowest root of $\Sigma_{\mathfrak{g}}$ with respect to the simple system $\Pi_{\mathfrak{g}}$. For each simple Lie algebra $\mathfrak{g}$ one can consult such data of Satake diagrams in \cite[Table 9 (pp. 312--316)]{ov2}.

\textbf{Step 3.} Let us fix a vertex $\gamma$ in $\mathrm{EDD}_{\Sigma}$. For example we take $\gamma=\lambda_{0}.$ We denote by $C$ the set of all vertices $\lambda$ ($\neq \gamma$) in $\mathrm{EDD}_{\Sigma}$ which are not connected to $\gamma$ in $\mathrm{EDD}_{\Sigma}$. In our situation, $C=\{  \lambda_{2},\lambda_{3} \}.$

\textbf{Step 4.} Let us define $\mathfrak{h} := [Z_{\mathfrak{g}}(A_{\gamma}),Z_{\mathfrak{g}}(A_{\gamma})].$ Then the Satake diagram $S_{\mathfrak{h}}$ is obtained by removing every white root $\alpha_{i}$ such that $r(\alpha_{i})\notin C$ from $S_{\mathfrak{g}}$. Thus in our example $S_{\mathfrak{h}}$ is of the following form

$$ 
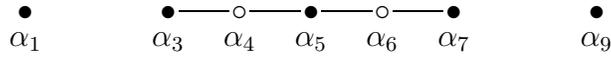
\xymatrix@1@R=2pt@!C=3pt{
{\bullet}  &  &{\bullet}  \ar@{-}[r]& {\circ} \ar@{-}[r] &{\bullet}  \ar@{-}[r]& {\circ}\ar@{-}[r] &{\bullet}  &  &{\bullet}\\ \alpha_1 & &\alpha_3&\alpha_4&\alpha_5&\alpha_6&\alpha_7&&\alpha_9
 }$$

Therefore, $\mathfrak{h}$ is isomorphic to $\mathfrak{su}(2)\oplus \mathfrak{su}(2)\oplus \mathfrak{su}^{*}(6)$.

\textbf{Step 5.} For a closed subgroup $H$ of $G$ with $\textrm{Lie} \ H = \mathfrak{h}$, the homogeneous space $G/H$ admits a proper $SL(2,\mathbb{R})$-action via $G$ by Theorem \ref{twggg1}.

As in Table \ref{tab1}, we give some examples of $(\mathfrak{g},\mathfrak{h})$ obtained by this algorithm.

\begin{center}
 \begin{table}[h]
 \centering
 {\footnotesize
 \begin{tabular}{| c | c |c| c |}
   \hline
   $\mathfrak{g}$ & $\mathfrak{h}$ & $\mathrm{EDD}_{\Sigma}$ & $\gamma$ \\
   \hline
   $\mathfrak{su}(p,q)$, $4<p<q$  & $\mathfrak{su}(p-1,q-1)$  & $\xymatrix@1@R=2pt@!C=3pt{
\lambda_{0} & \lambda_{1} & \lambda_{2} & \lambda_{p-1} & \lambda_{p}   \\
{\circledcirc}  \ar@{=>}[r]& {\circ}  \ar@{-}[r]& {\circ}   \ar@{.}[r]&{\circ}   \ar@{=>}[r]&{\circ} }$ & $\lambda_{0}$\\

   \hline
   $\mathfrak{su}(p,p)$, $4<p$ & $\mathfrak{su}(p-1,p-1)$  & $\xymatrix@1@R=2pt@!C=3pt{
\lambda_{0} & \lambda_{1} & \lambda_{2} & \lambda_{p-1} & \lambda_{p}   \\
{\circledcirc}  \ar@{=>}[r]& {\circ}  \ar@{-}[r]& {\circ}   \ar@{.}[r]&{\circ}   \ar@{<=}[r]&{\circ} }$ & $\lambda_{0}$\\

	\hline
   $\mathfrak{so}(p,q) $, $4<p$  & $\mathfrak{so}(p-2,q-2)$  & $  \quad\xymatrix@1@R=2pt@!C=3pt{   &\lambda_{0}&&&\\
  &{\circledcirc} \ar@{-}[dd]& & &\\
\lambda_{1} &   &\lambda_{3}&\lambda_{p-1} & \lambda_{p} \\
{\circ}  \ar@{-}[r]& {\circ} \ar@{-}[r]&{\circ} \ar@{.}[r]&{\circ}   \ar@{=>}[r]& {\circ}   \\
 & \lambda_{2} &&& } $ & $\lambda_{0}$\\
  \hline
     $\mathfrak{sp}(p,q)$, $4<p<q$  & $\mathfrak{sp}(p-1,q-1)$  & $\xymatrix@1@R=2pt@!C=3pt{
\lambda_{0} & \lambda_{1} & \lambda_{2} & \lambda_{p-1} & \lambda_{p}   \\
{\circledcirc}  \ar@{=>}[r]& {\circ}  \ar@{-}[r]& {\circ}   \ar@{.}[r]&{\circ}   \ar@{=>}[r]&{\circ} }$ & $\lambda_{0}$\\

   \hline
   $\mathfrak{sp}(p,p)$, $4<p$ & $\mathfrak{sp}(p-1,p-1)$  & $\xymatrix@1@R=2pt@!C=3pt{
\lambda_{0} & \lambda_{1} & \lambda_{2} & \lambda_{p-1} & \lambda_{p}   \\
{\circledcirc}  \ar@{=>}[r]& {\circ}  \ar@{-}[r]& {\circ}   \ar@{.}[r]&{\circ}   \ar@{<=}[r]&{\circ} }$ & $\lambda_{0}$\\

	   \hline
	  $\mathfrak{e}_{6}^{\textrm{II}}$ & $\mathfrak{sl}(3,\mathbb{R})$  & $\xymatrix@1@R=2pt@!C=3pt{
\lambda_{0} & \lambda_{1} & \lambda_{2} & \lambda_{3} & \lambda_{4}   \\
{\circledcirc}  \ar@{-}[r]& {\circ}  \ar@{-}[r]& {\circ}   \ar@{=>}[r]&{\circ}   \ar@{-}[r]&{\circ} }$ & $\lambda_{1}$\\
   \hline
	  $\mathfrak{e}_{7}^{\textrm{VI}}$ & $\mathfrak{so}(3,7)$  & $\xymatrix@1@R=2pt@!C=3pt{
\lambda_{0} & \lambda_{1} & \lambda_{2} & \lambda_{3} & \lambda_{4}   \\
{\circledcirc}  \ar@{-}[r]& {\circ}  \ar@{-}[r]& {\circ}   \ar@{=>}[r]&{\circ}   \ar@{-}[r]&{\circ} }$ & $\lambda_{0}$\\
   \hline
		  $\mathfrak{e}_{7}^{\textrm{VI}}$ & $\mathfrak{su}^{*}(6)$  & $\xymatrix@1@R=2pt@!C=3pt{
\lambda_{0} & \lambda_{1} & \lambda_{2} & \lambda_{3} & \lambda_{4}   \\
{\circledcirc}  \ar@{-}[r]& {\circ}  \ar@{-}[r]& {\circ}   \ar@{=>}[r]&{\circ}   \ar@{-}[r]&{\circ} }$ & $\lambda_{4}$\\
   \hline  
	$\mathfrak{e}_{7}^{\textrm{VII}}$ & $\mathfrak{so}(2,10)$  & $\xymatrix@1@R=2pt@!C=3pt{
\lambda_{0} & \lambda_{1} & \lambda_{2} & \lambda_{3}    \\
{\circledcirc}  \ar@{=>}[r]& {\circ}  \ar@{-}[r]& {\circ}   \ar@{<=}[r]&{\circ} }$ & $\lambda_{0}$\\
	\hline
 \end{tabular}
 }
\captionsetup{justification=centering}
 \caption{
 Examples of $(\mathfrak{g}, \mathfrak{h})$ obtained by the algorithm described in Section \ref{sub52}.
 }
 \label{tab1}
 \end{table}
\end{center}

\subsection{Homogeneous spaces with $\text{rank}_{\mathbb{R}}H=1$}\label{subsec:rank-one}

In this subsection, we study the cases where $\textrm{rank}_{\mathbb{R}} H=1$.

Let $G$ be a connected real semisimple linear Lie group
and $H$ a closed subgroup of $G$ which is reductive in $G$ (see Section \ref{sec:triples} for the notation).
We denote by $\mathfrak{g}$ and $\mathfrak{h}$ the Lie algebras of $G$ and $H$, respectively as before. 
Suppose that $\textrm{rank}_{\mathbb{R}}(\mathfrak{h})=1$ and the center of $H$ is compact. 
This means that $\theta$ is identity on the center of $\mathfrak{h}$.

\begin{proposition}\label{twgg}
In the setting above, $G/H$ admits a proper $SL(2,\mathbb{R})$-action via $G$ if and only if $\textrm{rank}_{\textrm{a-hyp}}\mathfrak{g} \geq 2$ $($see $\cite{bt}$ for the definition of the a-hyperbolic rank ``$\textrm{rank}_{\textrm{a-hyp}}\mathfrak{g}$'' of a semisimple Lie algebra $\mathfrak{g}$$)$.
In particular, the first condition for $G/H$ does not depend on the choice of $H (\subset G)$ in our setting.
\end{proposition} 

\begin{proof}
The proof is just a corollary to \cite{bt} and \cite{ok}. Indeed, the ``only if'' part follows from an observation that $\text{rank}_{\mathbb{R}}H=\text{rank}_{\textrm{a-hyp}}H=1$ and \cite[Theoorem 8]{bt}, the ``if part'' follows from the definition of the a-hyperbolic rank and \cite[Proposition 4.8 (i) and Corollary 5.5]{ok}.
\end{proof}

If $\rank_\R \mathfrak{g} = 1$, 
then the Calabi--Markus phenomenon occurs on $G/H$ (see \cite{kob2}).
Thus our interesting is in the cases where $\rank_\R \mathfrak{g} > 1$. 
By the computation of a-hyperbolic ranks of real simple Lie algebras in \cite{bt}, 
we have the following classification result.

\begin{corollary}\label{coo1}
In the setting above, assume that $G$ is real simple and $\text{rank}_{\mathbb{R}}(\mathfrak{g}) > 1$. 
Then $G/H$ admits a proper $SL(2,\mathbb{R})$-action via $G$ 
if and only if 
$\mathfrak{g}$ is not isomorphic to either 
$\mathfrak{sl}(3,\mathbb{R})$, $\mathfrak{su}^{*}(6)$ or $\mathfrak{e}^{\mathrm{IV}}_{6}$.
\end{corollary}

MB, PJ, AT: 

Department of Mathematics and Computer Science, University of Warmia and Mazury,

S\l\/oneczna 54, 10-710, Olsztyn, Poland

mabo@matman.uwm.edu.pl (MB), 

piojas@matman.uwm.edu.pl (PJ),

tralle@matman.uwm.edu.pl (AT).

TO:

Department of Mathematics, Hiroshima University,

1-3-1 Kagamiyama, Higashi-Hiroshima, 739-8526 Japan,

okudatak@hiroshima-u.ac.jp (TO).

\end{document}